%% file: ACGG-revised.tex
\documentclass[12pt]{article}
\usepackage{amsmath}
\usepackage{amssymb}
\usepackage[dvips]{graphicx}
\usepackage{epsfig}
\usepackage{vector}
\font\tencmmib=cmmib10 \skewchar\tencmmib '60
\newfam\cmmibfam
\textfont\cmmibfam=\tencmmib

\def\bbox{\quad\hbox{\vrule \vbox{\hrule \vskip2pt \hbox{\hskip2pt
\vbox{\hsize=1pt}\hskip2pt} \vskip2pt\hrule}\vrule}}
\def\lessim{\ \lower4pt\hbox{$
\buildrel{\displaystyle <}\over\sim$}\ }
\def\gessim{\ \lower4pt\hbox{$\buildrel{\displaystyle >}
\over\sim$}\ }

\def\la{{\Bigl\langle}}
\def\ra{{\Bigr\rangle}}

\def\qed{\hfill\break\rightline{$\bbox$}}
\newcommand{\e}{\mathbb{E}}

\newcommand{\Reals}{\mathbb{R}}
\newcommand{\Natural}{\mathbb{N}}

\newcommand{\vsi}{{\vec{\sigma}}}
\newcommand{\vrho}{{\vec{\rho}}}

\newtheorem{lemma}{Lemma}
\newtheorem{theorem}{Theorem}

\makeatletter
\@addtoreset{equation}{section}

\makeatother

\input invlat.tex

\begin{document}

\title{A unified  stability property in spin glasses.}
\author{Dmitry Panchenko\thanks{Department of Mathematics, Texas A\&M University, email: panchenk@math.tamu.edu. Partially supported by NSF grant.}}
\maketitle
\begin{abstract} 
Gibbs' measures in the Sherrington-Kirkpatrick type models satisfy two asymptotic stability properties,
the Aizenman-Contucci stochastic stability and the Ghirlanda-Guerra identities, which 
play a fundamental  role in our  current understanding of these models. In this paper we show that 
one can combine these two properties very naturally into one unified stability property.

\end{abstract}
\vspace{0.5cm}

\section{Introduction and main results.}

In the Sherrington-Kirkpatrick (SK) model \cite{SherK} or, more generally, in the mixed $p$-spin model
one considers a random process (Hamiltonian) $H_N(\vsi)$ indexed by spin configurations 
$\vsi\in\Sigma_N = \{-1,+1\}^N$ given by a linear combination $\sum_{p\geq 1} \beta_p H_{N,p}(\vsi)$
with some coefficients $\beta_p\geq 0$ of the independent Gaussian processes, called $p$-spin Hamiltonians,
\begin{equation}
H_{N,p}(\vsi)
=
\frac{1}{N^{(p-1)/2}}
\sum_{1\leq i_1,\ldots,i_p\leq N}g_{i_1,\ldots,i_p} \sigma_{i_1}\ldots\sigma_{i_p},
\label{HamSKp}
\end{equation}
where $(g_{i_1,\ldots,i_p})$ are standard Gaussian independent for all $p\geq 1$ and all $(i_1,\ldots,i_p).$ 
The Gibbs measure $G_N$ corresponding to the Hamiltonian $H_N$ is defined as  a random probability measure 
on $\Sigma_N$ given by 
\begin{equation}
G_N(\vsi) = \frac{1}{Z_N} \exp H_N(\vsi),
\label{Gibbs}
\end{equation} 
where the normalizing factor $Z_N$ is called the partition function. 
The Gibbs measure $G_N$ in (\ref{Gibbs}) is the central object of interest in spin glass models and the answers 
to many important questions follow from the conjectured properties of $G_N$. These properties can be expressed 
in terms of various functions of the sample $(\vsi^l)_{l\geq 1}$ from $G_N$, for example, the normalized Gram matrix 
$R = (R_{l,l'})_{l,l'\geq 1} = N^{-1}(\vsi^l \cdot \vsi^{l'})_{l,l'\geq 1}$. 
It is easy to see that knowing $R$ is equivalent to knowing $G_N$ up to orthogonal transformations since one can reconstruct 
$G_N$ from $R$ up to orthogonal transformations (as we shall see in the proof of Theorem \ref{Th1} below) and the information 
encoded in the distribution of $R$ turns out to be sufficient for most purposes in the setting of the SK model due to the fact that 
the Hamiltonian $H_N(\vsi)$ is a Gaussian process and its covariance $\e H_N(\vsi^1)H_N(\vsi^2)$ is the function of exactly 
the normalized scalar product $R_{1,2} = N^{-1}\vsi^1 \cdot \vsi^{2}$, called the overlap of $\vsi^1$ and $\vsi^2$.
Given any limiting distribution of the Gram matrix $R$ in the thermodynamic limit $N\to\infty$, 
one can use the Dovbysh-Sudakov representation (\cite{DS},\cite{PDS}) to define the asymptotic analogue 
of the Gibbs measure as a random probability measure $G$ on the unit ball of the Hilbert space $\ell_2$ (see \cite{AA}, \cite{ADS}, \cite{PGG}).
This means that in the limit the matrix $R$ is still generated as the Gram matrix of the sample from 
some random measure $G.$ For simplicity, let us assume that this asymptotic Gibbs measure $G$ is atomic, 
\begin{equation}
G=\sum_{l\geq 1} w_l \delta_{\xi_l},
\label{G}
\end{equation} 
with the weights arranged in non-increasing order, $w_1\geq w_2\geq \ldots,$ and let us
 denote by $Q=(\xi_l\cdot\xi_{l'})_{l,l'\geq 1}$ the matrix of scalar products of the points in the support of $G$. 
 Let $(\vsi^l)_{l\geq 1}$ again be an i.i.d. sample from this measure and let $R_{l,l'} = \vsi^l\cdot \vsi^{l'}$ be 
 the scalar product in $\ell_2$  of $\vsi^l$ and $\vsi^{l'}$. For any $n\geq 1$ and a function $f=f(\vsi^1,\ldots,\vsi^n)$ of $n$ 
configurations we will denote its average with respect to $G^{\otimes \infty}$ by
\begin{equation}
\la f \ra = \sum_{l_1,\ldots, l_n\geq 1} w_{l_1}\cdots w_{l_n} \,f(\xi_{l_1},\ldots,\xi_{l_n}).
\label{Av}
\end{equation}
We will denote by $\e$ the expectation with respect to the randomness of $G$. 

In the mixed $p$-spin models, one expects the asymptotic Gibbs measures (\ref{G}) to be described precisely 
by the Parisi ultrametric ansatz (see \cite{Parisi}, \cite{MPV}). 
 So far, most of the progress in the direction of proving structural results about $G$ was based on the idea that certain 
information about its geometry can be recovered from the asymptotic stability properties of the Gibbs measure $G_N$ under small 
perturbations of the Hamiltonian $H_N(\vsi)$.  In particular, two such stability properties in the setting of mixed $p$-spin model
are well known - the Aizenman-Contucci stochastic stability \cite{AC} and the Ghirlanda-Guerra identities \cite{GG}. 
They can be written down in terms of the asymptotic Gibbs measure $G$ as follows. 
\smallskip\\
1. {\it (Ghirlanda-Guerra identities)} 
Random measure $G$ is said to satisfy the Ghirlanda-Guerra identities \cite{GG} if for any $n\geq 2,$ any bounded measurable function $f$ 
of the overlaps  $(R_{l,l'})_{l,l'\leq n}$ and any integer $p\geq 1$ we have
\begin{equation}
\e \la f R_{1,n+1}^p\ra = \frac{1}{n}\, \e\la f \ra \e\la R_{1,2}^p\ra + \frac{1}{n}\sum_{l=2}^{n}\e\la f R_{1,l}^p\ra.
\label{GG}
\end{equation}
These constraints on the distribution of $R$ look very mysterious but they arise from a very natural and very general principle 
of the concentration of the Hamiltonian $H_N(\vsi)$ (see \cite{GG}). 
\smallskip \\  
2. {\it (Aizenman-Contucci stochastic stability)}
Given integer $p\geq 1,$ let $(g_p(\xi_l))_{l\geq 1}$ be a Gaussian sequence conditionally on $G$ indexed
by the points $(\xi_l)_{l\geq 1}$ with covariance
\begin{equation}
\mbox{\rm Cov}\bigl(g_p(\xi_l),g_p(\xi_{l'})\bigr) = (\xi_l\cdot\xi_{l'})^p.
\label{Covp}
\end{equation}
Given $t\in\Reals,$ consider a new sequence of weights
\begin{equation}
w_l^t=\frac{w_l e^{t g_p(\xi_l)}}{\sum_{j\geq 1} w_j e^{t g_p(\xi_j)}}
\label{density}
\end{equation}
defined by a random change of density proportional to $e^{tg_p(\xi_l)}.$ Let $(w_l^\pi)$ be the weights $(w_l^t)$ 
arranged in the non-increasing order and let $\pi:\Natural\to \Natural$ be the permutation keeping track of where
each index came from, $w_{l}^\pi = w_{\pi(l)}^t.$ Let us define by
\begin{equation}
G^\pi =\sum_{l\geq 1} w_l^{\pi} \delta_{\xi_{\pi(l)}}
\,\,\mbox{ and }\,\,
Q^\pi= \bigl(\xi_{\pi(l)}\cdot  \xi_{\pi(l')}\bigr)_{l,l'\geq 1}
\label{change}
\end{equation}
the probability measure $G$ after the change of density proportional to $e^{tg_p(\xi_l)}$
and the matrix $Q$ rearranged according to the reordering of weights. Measure $G$ is said to satisfy
the Aizenman-Contucci stochastic stability \cite{AC} if for any $p\geq 1$ and $t\in \Reals,$
\begin{equation}
\bigl(
(w_l^\pi)_{l\geq 1}, Q^\pi
\bigr)
\stackrel{d}{=}
\bigl(
(w_l)_{l\geq 1}, Q
\bigr)
\label{AC}
\end{equation}
where equality in distribution is in the sense of finite dimensional distributions  of these arrays.
This property represents the invariance of the distribution of measure $G$ up to orthogonal
transformations under the random changes of density (\ref{density}) and it arises from the continuity
of the Gibbs measure $G_N$ under small changes in the inverse temperature parameters $\beta_p$
in the Hamiltonian $H_N(\vsi)$ (see \cite{AC}).
\qed 

Originally, the proofs of the Ghirlanda-Guerra identities (\ref{GG}) in \cite{GG} and the Aizenman-Contucci stochastic 
stability (\ref{AC}) in \cite{AC}, \cite{CGSS} obtained these results for each $p\geq 1$ on average of the parameter $\beta_p$ over 
any non-trivial interval. In a closely related formulation, one can always perturb the parameters $(\beta_p)$ slightly 
(for example, one can find $(\beta_{N,p})$ such that $|\beta_{N,p} - \beta_p |\leq 2^{-p}N^{-1/16}$) so that the sequence
of Gibbs measures $G_N$ corresponding to these slightly perturbed parameters satisfies the above properties in
the limit (see \cite{Tal-New}, \cite{SG2}). Both of these formulations hold for mixed $p$-spin models with arbitrary subset
of $p$-spin terms (\ref{HamSKp}) present in the model, i.e. for which $\beta_p\not = 0$. However, if the model contains terms
for $p=1$ and all even $p\geq 2$ then it was proved in \cite{PGGmixed} that the Ghirlanda-Guerra identities (\ref{GG}) hold in the
strong sense without perturbation of parameters for all integer $p\geq 1$. The proof was based on the validity of the Parisi formula 
for the free energy proved by M. Talagrand in \cite{T-P} following the discovery of the replica symmetry breaking bound by
F. Guerra in \cite{Guerra}, on the differentiability properties of the Parisi formula (see \cite{PM}, \cite{ECPs}) and on the
positivity of the overlap (see \cite{SG2}). A similar strong version of the Aizenman-Contucci stochastic stability  was proved 
in \cite{ACh}.

As we mentioned above, the importance of the stability properties (\ref{GG}) and (\ref{AC}) in the SK model comes from their many applications 
(see e.g. \cite{AGG}, \cite{AA},  \cite{ACh}, \cite{PGG}, \cite{Posit}, \cite{PGGmixed}, \cite{PDI},  \cite{PGGsimple}, \cite{PGGnew}, 
\cite{PT}, \cite{RA}, \cite{Tal-New},  \cite{SG}, \cite{SG2}). 
Of course, the ultimate goal would be to prove the Parisi ultrametricity conjecture which states that the support 
of the asymptotic Gibbs measure $G$ must be ultrametric in $\ell_2$ with probability one, which would allow us 
to identify $G$ with the Ruelle Probability Cascades in \cite{Ruelle}. At the moment, the results that come closest to proving 
this conjecture are based either on the Aizenman-Contucci stochastic stability or on the Ghirlanda-Guerra identities.
First such result was proved by L.-P. Arguin and M. Aizenman in \cite{AA} using the Aizenman-Contucci stochastic stability
under a technical assumption that the scalar products $\xi_l \cdot \xi_{l'}$ of points in the support of measure $G$ take finitely 
many non-random values. Following their work, a similar result was proved by the author in \cite{PGG} (see also \cite{PGGsimple}
for a recent elementary proof)  and by M. Talagrand in \cite{Tal-New} using the Ghirlanda-Guerra identities instead. 
Some modest progress toward the general case was made in \cite{PGGnew} but the conjecture still remains open. 
Once this conjecture is proved,  the Parisi formula for the free energy will easily follow using the Aizenman-Sims-Starr scheme 
developed in \cite{AS2} which would naturally complete the mathematical justification of the Parisi ansatz in the SK model. 
It is worth mentioning that stability properties of the Gibbs measure under
small perturbations of the Hamiltonian play very important role in other spin glass models as well (see e.g. \cite{Pspins}).

In the main result of this paper we will show that one can combine the Ghirlanda-Guerra identities (\ref{GG}) and the Aizenman-Contucci
stochastic stability (\ref{AC}) into a joint stability property as follows. It is known (Theorem 2 in \cite{PGG}) that if the measure $G$ satisfies the Ghirlanda-Guerra identities
and if $q^*$ is the supremum of the support of the distribution of the overlap $R_{1,2}$ under $\e G^{\otimes 2}$
then with probability one $G$ is concentrated on the sphere of radius $\sqrt{q^*}$. Let
\begin{equation}
b_p = (q^*)^p - \e\la R_{1,2}^p\ra.
\label{bp}
\end{equation}
Then the following holds.
\begin{theorem}\label{Th1}
A random measure $G$ satisfies the Ghirlanda-Guerra identities (\ref{GG}) and the Aizenman-Contucci stochastic stability (\ref{AC}) 
if and only if  it is concentrated on the sphere of constant radius $\sqrt{q^*}$  with probability one and for any $p\geq 1$ and $t\in\Reals,$
\begin{equation}
\Bigl(
\bigl(w_l^\pi\bigr)_{l\geq 1}, \bigl(g_p(\xi_{\pi(l)}) - b_p t\bigr)_{l\geq 1}, Q^\pi
\Bigr)
\stackrel{d}{=}
\Bigr(
\bigl(w_l\bigr)_{l\geq 1}, \bigl(g_p(\xi_l)\bigr)_{l\geq 1}, Q
\Bigr).
\label{main}
\end{equation}
where equality in distribution is in the sense of finite dimensional distributions.
\end{theorem}
One can see that the Ghirlanda-Guerra identities are now replaced by the statement that the Gaussian field $(g_p(\xi_l))$ after
permutation $\pi$ corresponding to the reordering of weights in (\ref{density}) will only differ by a constant shift $b_p t$ in distribution.
In the language of competing particle systems (\cite{RA} and \cite{AA}), (\ref{main}) means that the past increments of the dynamics 
after re-centering have the same law as the future or forward increments.
The  stability property (\ref{main}) is well-known for the ultrametric Ruelle Probability Cascades (see Theorem 4.2 in \cite{Bolthausen}  
or Theorem 15.2.1 in \cite{SG2}) and, of course, the big question is whether it holds only for these measures and whether (\ref{main})
implies ultrametricity.  

The Ghirlanda-Guerra identities do not require the random measure $G$ to be discrete and, in fact, the Aizenman-Contucci 
stochastic stability can be formulated not only for discrete measures as well. We will mention this more 
general formulation in the next section. However, we prefer to state our main result in the setting of discrete measures 
since it allows for a particularly attractive formulation (\ref{main}) in the spirit of competing particle systems, as in \cite{RA} and \cite{AA}.
Moreover, from the point of view of studying structural properties of such measures one can without loss of
generality start with discrete measures since it is easy to show that sampling an i.i.d. sequence of points 
from the original measure and assigning them new independent weights from the Poisson-Dirichlet distribution 
creates a discrete measure which still satisfies both properties. On the other hand, almost any geometric
property of the original measure will be encoded into a countable i.i.d. sample and, therefore, this new discrete
measure. 

The unified stability property (\ref{main}) inspired a new representation of the Ghirlanda-Guerra identities in \cite{PGGnew}
which yielded some interesting applications.  For another recent stability property that reproduces the Ghirlanda-Guerra 
identities on average see \cite{CGSSCG}.
\medskip\\
\textbf{Acknowledgment.}  The author would like to thank the referees for making many important suggestions that 
helped improve the paper.

\section{Proof of Theorem \ref{Th1}.}

Let $(\vrho^l)_{l\geq 1}$ be an i.i.d. sequence from measure $G^\pi$ defined in (\ref{change}) and 
denote by $S_{l,l'} = \vrho^l\cdot \vrho^{l'}$ the overlap of $\vrho^l$ and $\vrho^{l'}$. Analogously to (\ref{Av}),
for any $n\geq 1$ and a function $f=f(\vrho^1,\ldots,\vrho^n)$ of $n$ configurations we will denote its average 
with respect to $(G^\pi)^{\otimes \infty}$ by
\begin{equation}
\la f \ra_\pi = \sum_{l_1,\ldots, l_n\geq 1} w_{l_1}^\pi\cdots w_{l_n}^\pi  \,f(\xi_{\pi(l_1)},\ldots,\xi_{\pi(l_n)}).
\label{Av2}
\end{equation}
We now will denote by $\e$ the expectation with respect to the randomness of $G$ and the Gaussian sequence $(g_p).$
Let us first make a simple observation that equality of finite dimensional distributions in (\ref{AC}) and (\ref{main}) implies
equality of averages with respect to the random measures in the following sense.
\begin{lemma}\label{Lem1}
If (\ref{main}) holds then for any $k\geq 1,$ any bounded measurable function $f$ of the overlaps on $k$ replicas and any integers
$n_1,\ldots, n_k \geq 0,$
\begin{equation}
\e\Bigl\la \prod_{l\leq k }\bigl(g_p(\vrho^l)-b_p t\bigr)^{n_l} f\bigl((S_{l,l'})_{l,l'\leq k}\bigr)\Bigr\ra_\pi
=
\e\Bigl\la \prod_{l\leq k } g_p(\vsi^l)^{n_l} f\bigl((R_{l,l'})_{l,l'\leq k}\bigr)\Bigr\ra
\label{Aveq}
\end{equation}
Under (\ref{AC}), this holds with all $n_l=0.$
\end{lemma}
\textbf{Remark.} One can consider (\ref{Aveq}) with all $n_l=0$ as the definition of the Aizenman-Contucci stochastic
stability for non-atomic measures in which case $(g_p(\xi))$ is the Gaussian field with covariance (\ref{Covp}).
Moreover, in this case (\ref{Aveq}) should be considered as the analogue of (\ref{main}).\medskip\\
\textbf{Proof.} This is obvious by separating the sum in (\ref{Av}) and (\ref{Av2}) into finitely many terms corresponding
to the largest weights and the remaining small weights. For example,
\begin{eqnarray*}
\e\Bigl\la \prod_{l\leq k } g_p(\vsi^l)^{n_l} f\bigl((R_{l,l'})_{l,l'\leq k}\bigr)\Bigr\ra
&=&
\e \sum_{j_1,\ldots, j_k\geq 1} w_{j_1}\cdots w_{j_n}
 \prod_{l\leq k } g_p(\xi_{j_l})^{n_l} f\bigl((\xi_{j_l}\cdot \xi_{j_{l'}})_{l,l'\leq k}\bigr)
 \\
&=&
\e \sum_{j_1,\ldots, j_k \leq N} w_{j_1}\cdots w_{j_n}
 \prod_{l\leq k } g_p(\xi_{j_l})^{n_l} f\bigl((\xi_{j_l}\cdot \xi_{j_{l'}})_{l,l'\leq k}\bigr) 
 +{\cal R}_N,
\end{eqnarray*}
where the remainder ${\cal R}_N$ consists of the terms with at least one index $j_1,\ldots, j_k> N$.
The left hand side of (\ref{Aveq}) can be similarly broken into two sums. The finite sums are equal because they
involve only finitely many elements of the arrays (\ref{main}) which are equal in distribution by assumption.
Thus, we only need to show that ${\cal R}_N$ becomes small for large $N$. First taking expectation in the Gaussian  
random variables $(g_p(\xi_l))$ conditionally on $(w_l)$ and $Q$ and using that
$$
\e \Bigl( \prod_{l\leq k } |  g_p(\xi_{j_l}) | ^{n_l} \,\Bigr|\, (w_l), Q\Bigr) \leq L(n_1,\ldots,n_k)
$$
we get that
$$
|{\cal R}_N| \leq L(n_1,\ldots,n_k) \|f\|_{\infty}\, \e \sum_{(j_1,\ldots, j_k \leq N)^c} w_{j_1}\cdots w_{j_n}
\leq
Lk\, \e\sum_{j> N} w_j
$$
which goes to zero as $N\to\infty.$ The remainder ${\cal R}_N^\pi$ for the left hand side of (\ref{Aveq}) is controlled
by exactly the same bound because, by (\ref{main}), the corresponding terms in ${\cal R}_N^\pi$ and ${\cal R}_N$
are equal in distribution.
\qed

The ``if" part of the Theorem \ref{Th1} is easy since assuming (\ref{main}) we only need to prove (\ref{GG}) 
and this follows from integration by parts of (\ref{Aveq}) with $n_1=1,n_2=\ldots=n_k=0.$ In this case the right
hand side is zero by averaging $g_p(\vsi^1)$ first and the left hand side is 
\begin{eqnarray*}
\e\Bigl\la (g_p(\vrho^1)-b_p t\bigr) f\bigl((S_{l,l'})_{l,l'\leq k}\bigr)\Bigr\ra_\pi
&=&
t \e\Bigl\la 
\bigl(
\sum_{l=1}^k S_{1,l}^p - b_p - kS_{1,k+1}^p
\bigr)
f\bigl((S_{l,l'})_{l,l'\leq k}\bigr)\Bigr\ra_\pi
\\
&=&
t \e\Bigl\la 
\bigl(
\sum_{l=1}^k R_{1,l}^p -b_p - k R_{1,k+1}^p
\bigr)
f\bigl((R_{l,l'})_{l,l'\leq k}\bigr)\Bigr\ra
\\
&=&
t \e\Bigl\la 
\bigl(
\sum_{l=2}^k R_{1,l}^p +\e\la R_{1,2}^p\ra - kR_{1,k+1}^p
\bigr)
f\bigl((R_{l,l'})_{l,l'\leq k}\bigr)\Bigr\ra
\end{eqnarray*}
where in the second line we used (\ref{AC}) part of (\ref{main}) and Lemma \ref{Lem1}, and in the third line
we used  (\ref{bp}) and the fact that  $\xi_l\cdot\xi_l = q^*$. The fact that the last sum is zero is exactly (\ref{GG}).
To prove the "only if" part we need the following key lemma.
\begin{lemma}\label{Lem2}
If (\ref{GG}) and (\ref{AC}) hold then (\ref{Aveq}) holds.
\end{lemma}
\textbf{Proof.} The proof is by induction on $N=n_1+\ldots+n_k$. When $N=0,$ (\ref{Aveq}) is the consequence
of (\ref{AC}) by Lemma \ref{Lem1}. Suppose (\ref{Aveq}) holds for all $k\geq 1$, all $f$  and for all $N\leq N_0.$
Clearly, we only need to prove the case of powers $n_1+1, n_2,\ldots, n_k.$ Writing 
$$
g_p(\vsi^1)^{n_1+1} = g_p(\vsi^1)\, g_p(\vsi^1)^{n_1} 
$$ 
and using Gaussian integration by parts for $g_p(\vsi^1)$ we can rewrite the right hand side of (\ref{Aveq}) 
with $n_1+1$ instead of $n_1$ as
\begin{equation}
\sum_{l\leq k} n_l 
\e\Bigl\la g_p(\vsi^1)^{n_1}\ldots g_p(\vsi^l)^{n_l-1}\ldots g_p(\vsi^k)^{n_k} 
R_{1,l}^p\, f\bigl((R_{l,l'})_{l,l'\leq k}\bigr)\Bigr\ra.
\label{rhs}
\end{equation}
Again, writing 
$$
(g_p(\vrho^1)-b_p t)^{n_1+1} = (g_p(\vrho^1)-b_p t) (g_p(\vrho^1)-b_p t)^{n_1} 
$$ 
and using Gaussian integration by parts for $g_p(\vrho^1)$ we can rewrite the left hand side of (\ref{Aveq})
with $n_1+1$ instead of $n_1$ as $\rm I + II$ where $\rm I$ is given by
\begin{equation}
\sum_{l\leq k} n_l 
\e\Bigl\la \bigl(g_p(\vrho^1)-b_p t\bigr)^{n_1}\ldots \bigl(g_p(\vrho^l)-b_p t\bigr)^{n_l-1}\ldots 
\bigl(g_p(\vrho^k)-b_p t\bigr)^{n_k} 
S_{1,l}^p\, f\bigl((S_{l,l'})_{l,l'\leq k}\bigr)\Bigr\ra_\pi
\label{lhsI}
\end{equation}
and $\rm II$ is given by
\begin{equation}
t \, \e\Bigl\la \prod_{l\leq k }\bigl(g_p(\vrho^l)-b_p t\bigr)^{n_l} 
\Bigl(
\sum_{l\leq k} S_{1,l}^p - b_p - k S_{1,k+1}^p
\Bigr)
f\bigl((S_{l,l'})_{l,l'\leq k}\bigr)\Bigr\ra_\pi.
\label{lhsII}
\end{equation}
By induction hypothesis, (\ref{lhsI}) is equal to (\ref{rhs}) and (\ref{lhsII}) is equal to
\begin{equation}
t \, \e\Bigl\la \prod_{l\leq k } g_p(\vsi^l)^{n_l} 
\Bigl(
\sum_{l\leq k} R_{1,l}^p - b_p - k R_{1,k+1}^p
\Bigr)
f\bigl((R_{l,l'})_{l,l'\leq k}\bigr)\Bigr\ra.
\label{lhsII2}
\end{equation}
Since $\la\cdot\ra$ does not depend on the Gaussian sequence $(g_p(\xi_l))$ we can take expectation
$\e_g$ with respect to the randomness of this sequence conditionally on $G$ first and notice that 
$$
\e_g \prod_{l\leq k } g_p(\vsi^l)^{n_l} = f'\bigl((R_{l,l'})_{l,l'\leq k}\bigr)
$$
for some function $f'$ of the overlaps of $k$ configurations $\vsi^1,\ldots,\vsi^k.$ Therefore, (\ref{lhsII2}) equals
\begin{eqnarray}
&&
t \, \e\Bigl\la
\Bigl(
\sum_{l=1}^k R_{1,l}^p - b_p - k R_{1,k+1}^p
\Bigr)
(ff')\bigl((R_{l,l'})_{l,l'\leq k}\bigr)\Bigr\ra
\\
&&=
t \, \e\Bigl\la
\Bigl(
\sum_{l=2}^k R_{1,l}^p + \e\la R_{1,2}^p\ra - k R_{1,k+1}^p
\Bigr)
(ff')\bigl((R_{l,l'})_{l,l'\leq k}\bigr)\Bigr\ra =0,
\nonumber
\end{eqnarray}
where in the first equality we again used  (\ref{bp}) and the fact that  $\xi_l\cdot\xi_l = q^*$ and
the second equality is by the Ghirlanda-Guerra identities (\ref{GG}). This finishes the proof.
\qed

The equality of joint moments (\ref{Aveq}) proved in Lemma \ref{Lem2} implies the following.
\begin{lemma}
If (\ref{GG}) and (\ref{AC}) hold then 
\begin{equation}
\Bigl(
\bigl( g_p(\vrho^l) - b_p t \bigr)_{l\geq 1}, (S_{l,l'})_{l,l'\geq 1}
\Bigr)
\stackrel{d}{=}
\Bigl(
\bigl( g_p(\vsi^l) \bigr)_{l\geq 1}, (R_{l,l'})_{l,l'\geq 1}
\Bigr)
\label{samples}
\end{equation}
where equality in distribution is in the sense of finite dimensional distributions.
\end{lemma}
\textbf{Remark.}
Let us recall that the i.i.d. sequences $(\vsi^l)$ and $(\vrho^l)$ are sampled from $G$ and $G^\pi$ correspondingly
and, therefore, the distributions of the right-hand side and left-hand side in (\ref{samples}) are under $\e G^{\otimes \infty}$
and $\e (G^\pi)^{\otimes \infty}$.\medskip\\
\textbf{Proof.}
By choosing $f$ to be monomials, (\ref{Aveq}) gives the equality of joint moments of 
the corresponding elements of the two arrays in (\ref{samples}). In our case the joint moments uniquely
determine joint distributions, for example, by the main result in \cite{Petersen} which states that we only 
need to ensure the uniqueness of one dimensional marginals and the fact that the one dimensional
marginals are either bounded or Gaussian.
\qed\medskip\\
\textbf{Proof of Theorem \ref{Th1}.}
Finally, we will show that (\ref{samples}) implies (\ref{main}).  The procedure is very similar to the one at the end of Theorem $4$ 
in \cite{PGG} or a more general argument in Lemma 4 in \cite{PDS}. First of all, by the well-known result of Talagrand
(Section 1.2 in \cite{SG}), the Ghirlanda-Guerra identities imply that the weights $(w_l)$ must have a Poisson-Dirichlet
distribution $PD(m)$ where $m$ is determined by
$$
\e\la I(R_{1,2}=q^*)\ra =\e\sum_{l\geq 1} w_l^2 = 1-m.
$$
This means that if $(u_l)$ is a Poisson point process on $(0,\infty)$ with intensity
measure $x^{-m-1}dx$ then $w_l = u_l/\sum_{j\geq 1} u_j.$ In particular, all the weights are different
with probability one. This point is not crucial but it makes for an easier argument.
The reason why (\ref{samples}) implies (\ref{main}) is because one can easily reconstruct the arrays
in (\ref{main}) from the arrays (\ref{samples}) using that $(\vsi^l)$ is an i.i.d. sample from $(\xi_l)$
according to weights $(w_l)$ and $(\vrho^l)$ is an i.i.d. sample from $(\xi_{\pi(l)})$ according to weights $(w_l^\pi).$
The key observation here is that given arrays (\ref{samples}) we know exactly when $\vsi^l = \vsi^{l'}$ and
$\vrho^l=\vrho^{l'}$ since this is equivalent to $R_{l,l'}=q^*$ and $S_{l,l'}=q^*$. Therefore, given $N\geq 1$ and
$$
\bigl(
\bigl( g_p(\vsi^l) \bigr)_{l\leq N}, (R_{l,l'})_{l,l'\leq N}
\bigr)
$$
we can partition the set $\{1,\ldots,N\}$ according to the equivalence relation $l\sim l'$ defined by $R_{l,l'}=q^*,$
let the sequence of weights $(w_l^N)_{l\geq 0}$ be the proportions of the sets in this partition arranged
in non-increasing order and extended by zeros and, given any integer $j$ in the element of the partition 
corresponding to the weight $w_l^N$, define $\xi_l^N = \vsi^j.$ We let $Q^N=(\xi_l^N\cdot\xi_{l'}^N)_{l,l'\geq 1}$.
The elements of $(\xi_l^N)$ and $Q^N$ with indices corresponding to zero weights $w_l^N$ can be set to some fixed 
values, and we break ties between $w_l^N$ by any pre-determined rule. Similarly, given
$$
\bigl(
\bigl( g_p(\vrho^l) - b_p t \bigr)_{ l\leq N}, (S_{l,l'})_{ l,l'\leq N}
\bigr)
$$
we can construct sequences $(\tilde{w}_l^N), (\tilde{\xi}^N_l)$ and $\tilde{Q}^N = (\tilde{\xi}_l^N\cdot \tilde{\xi}_{l'}^N)$.
Equation (\ref{samples}) implies that for any fixed $k\geq 1,$
$$
\Bigr(
\bigl(\tilde{w}^N_l\bigr)_{l\leq k}, \bigl( g_p(\tilde{\xi}_l^N) -b_p t \bigr)_{ l\leq k}, (\tilde{q}^N_{l,l'})_{ l,l'\leq k}
\Bigr)
\stackrel{d}{=}
\Bigr(
\bigl(w^N_l\bigr)_{ l\leq k}, \bigl(g_p(\xi^N_l) \bigr)_{ l\leq k}, (q^N_{l,l'})_{ l,l'\leq k}
\Bigr).
$$
It remains to observe that the right hand side converges
\begin{equation}
\Bigr(
\bigl(w^N_l\bigr)_{ l\leq k}, \bigl(g_p(\xi^N_l)\bigr)_{ l\leq k}, (q^N_{l,l'})_{ l,l'\leq k}
\Bigr)
\to
\Bigr(
\bigl(w_l\bigr)_{ l\leq k}, \bigl(g_p(\xi_l)\bigr)_{ l\leq k}, (q_{l,l'})_{ l,l'\leq k}
\Bigr)
\label{consist}
\end{equation}
almost surely and, similarly, the left hand side converges a.s. to the corresponding
array from the left hand side of (\ref{main}). To prove (\ref{consist}), we notice that by construction
$$
G^N:=\sum_{l\geq 1} w_l^N \delta_{\xi_l^N} = \frac{1}{N}\sum_{i\leq N} \delta_{\sigma_i}
$$
is the empirical measure based on the sample $\vsi^1,\ldots,\vsi^N$ from the measure 
$G=\sum_{l\geq 1}w_l \delta_{\xi_l}.$ By the strong law of large number for empirical measures
(e.g. Theorem 11.4.1 in \cite{Dudley}), the laws $G^N\to G$ almost surely and since the Poisson-Dirichlet
weights $(w_l)$ are all different a.s., the largest $k$ weights must converge
$(w^N_l)_{ l\leq k} \to (w_l)_{ l\leq k}$  almost surely and for large enough $N$ we must have 
$(\xi_l^N)_{l\leq k} = (\xi_l)_{l\leq k}$ and, thus, (\ref{consist}) holds.
\qed

\end{document}

%% file: invlat.tex
\font\tencmmib=cmmib10 \skewchar\tencmmib '60
\newfam\cmmibfam
\textfont\cmmibfam=\tencmmib

\def\bbox{\quad\hbox{\vrule \vbox{\hrule \vskip2pt \hbox{\hskip2pt
\vbox{\hsize=1pt}\hskip2pt} \vskip2pt\hrule}\vrule}}
\def\lessim{\ \lower4pt\hbox{$
\buildrel{\displaystyle <}\over\sim$}\ }
\def\gessim{\ \lower4pt\hbox{$\buildrel{\displaystyle >}
\over\sim$}\ }

\def\go0{\to 0}

\def\la{\langle}

\def\leftitem#1{\item{\hbox to\parindent{\enspace#1\hfill}}}

\def\qed{\hfill\break\rightline{$\bbox$}}

\def\ra{\rangle}

\def\sg{\sigma}

\def\sg2{\sigma^2}

\def\__{_{\infty}}